\documentclass[12pt,reqno]{amsart}
\theoremstyle{plain}    
\newtheorem{thm}{Theorem}[section]
\numberwithin{figure}{section} %% Comment out for sequentially-numbered
\theoremstyle{plain}

\newtheorem{cor}[thm]{Corollary} %%Delete [thm] to re-start numbering
\newtheorem{prop}[thm]{Proposition}
\newtheorem{numbered}[thm]{}

\newtheorem{remark}[thm]{Remark}

\usepackage{amscd,amssymb,comment,euscript}
\newcommand\Ac{{\mathcal{A}}}

\newcommand\AEu{{\EuScript A}}

\newcommand\clspan{{\overline{\mathrm{span}}\,}}

\newcommand\Cpx{{\mathbf C}}

\newcommand\Fb{{\mathbf F}}

\newcommand\fdim{\text{\rm fdim}\,}

\newcommand\GEu{{\EuScript G}}

\newcommand\HEu{{\EuScript H}}                   % requires package euscript

\newcommand\id{{\operatorname{id}}}

\newcommand\Ints{{\mathbf Z}}

\newcommand\Lambdao{{\Lambda\oup}}

\newcommand\lspan{\mathrm{span}\,}

\newcommand\Mcal{{\mathcal{M}}} %  needed to be renamed, because was ``\Mc already defined''
\newcommand\MEu{{\EuScript M}}                   % requires package euscript

\newcommand\MEuh{{\hat\MEu}}                   % requires package euscript

\newcommand\Mh{{\hat M}}

\newcommand\Nats{{\mathbf N}}

\newcommand\NEu{{\EuScript N}}                   % requires package euscript

\newcommand\Nh{{\hat N}}

\newcommand\oup{^{\mathrm o}}

\newcommand\Pc{{\mathcal{P}}}

\newcommand\PEu{{\EuScript P}}

\newcommand\Phit{{\widetilde\Phi}}

\newcommand\Psit{{\widetilde\Psi}}

\newcommand\Qc{{\mathcal{Q}}}

\newcommand\QEu{{\EuScript Q}}

\newcommand\Qt{{\widetilde Q}}

\newcommand\Reals{{\mathbf R}}

\newcommand\restrict{{\upharpoonright}}

\newcommand\REu{{\EuScript R}}                   % requires package euscript

\newcommand\staropwl{
  \operatornamewithlimits{\raisebox{-0.5ex}[1.5ex][0ex]{\rm*}}}

\newcommand\Tr{{\mathrm{Tr}}}

\begin{document}

\pagestyle{myheadings}

\title{Subfactors of free products of rescalings of a II$_1$--factor}
 
\author{Ken Dykema}

\thanks{Supported in part by
NSF grant DMS--0070558.
The author thanks the Mathematical Sciences Research Institute,
where part of this work was done.
Research at MSRI is supported in part by NSF grant DMS--9701755.}

\address{\hskip-\parindent
Department of Mathematics\\
Texas A\&M University\\
College Station TX 77843--3368, USA}
\email{Ken.Dykema@math.tamu.edu}

\begin{abstract}
Let $Q$ be any II$_1$--factor.
It is shown that any standard lattice $\GEu$ can be realized as the standard invariant of a free product of
(several) rescalings of $Q$.
In particular, if $Q$ has fundamental group equal to the positive reals and if $P$ is the free product
of infinitely many copies of $Q$,
then $P$ has subfactors giving rise to all possible standard invariants.
Similarly, given a II$_1$--subfactor $N\subset M$, it is shown there are subfactors $\Nh\subset\Mh$
having the same standard invariant as $N\subset M$ but where $\Mh$, respectively $\Nh$, is the free product
of $M$, respectively $N$, with rescalings of $Q$.
\end{abstract}

\date{25 February 2002}

\maketitle

\markboth{\tiny Subfactors in free products}{\tiny Subfactors in free products}
 
\section{Introduction}

The systems of higher relative commutants of finite index subfactors of II$_1$--factors were classified by
S.\ Popa~\cite{P}, when he proved certain conditions on a lattice of finite dimensional commuting squares (together with certain other data
such as Markov traces) to be
equivalent to it arising as the higher relative commutants of a subfactor.
Such an abstract lattice is called a {\em standard lattice}.
The question of whether there is a single II$_1$--factor $\Mcal$ whose subfactors give rise to all standard lattices
was answered in the affirmative by Popa and Shlyakhtenko~\cite{PS}, where they proved this universal property
for $\Mcal=L(\Fb_\infty)$,
the von Neumann algebra of the free group on infinitely many generators.
Their proof has two main components.
First, given a standard lattice $\GEu$, they construct a non--degenerate commuting square
\begin{equation}
\label{eq:cs}
\begin{matrix}
B    & \subset & A    \\
\cup &         & \cup \\
D    & \subset & C
\end{matrix}
\end{equation}
of type~I von Neumann algebras having atomic centers, which encodes the standard lattice.
Then, for any II$_1$--factor $Q$, they consider the resulting inclusion $N\subset M$ of amalgamated free product von Neumann algebras
\begin{align}
\MEu&=(Q\otimes B,\tau_Q\otimes\id_B)*_B(A,E^A_B) \label{eq:MNamalgfp} \\
\NEu&=(Q\otimes D,\tau_Q\otimes\id_D)*_D(C,E^C_D)\;, \notag
\end{align}
where $\tau_Q$ is the tracial state on $Q$ and $E^A_B$ and $E^C_D$ are the conditional expectations from the commuting square~\eqref{eq:cs}.
By relating it to a construction from~\cite{P}, they show that $\NEu\subset\MEu$ is an inculsion of II$_\infty$ factors which is an infinite
amplification of an inclusion $N\subset M$ of II$_1$--factors whose system of higher
relative commutants is the original standard lattice $\GEu$.
The second main part of their proof is to show that if $Q$ is taken to be $L(\Fb_\infty)$, then both $N$ and $M$
are isomorphic to $L(\Fb_\infty)$.

In this paper, we identify $N$ and $M$ for general $Q$ as a free scaled product~\cite{D} of rescalings of $Q$.
This allows us to generalize Popa and Shlyakhtenko's results.
For example, (Theorem~\ref{thm:subfin}), we show that for every standard lattice $\GEu$
of finite depth and for every II$_1$--factor $Q$ there is a subfactor $N\subseteq M$ whose standard lattice is $\GEu$
and where
\begin{align*}
M&=Q_{s_1}*\cdots*Q_{s_m} \\
N&=Q_{t_1}*\cdots*Q_{t_n}*L(\Fb_b)
\end{align*}
for some  $m,n\in\Nats$, positive real numbers $s_1,\ldots,s_m$ and $t_1,\ldots,t_n$
and $b>1-n$ depending only on $\GEu$.
(Here we are using the notation of~\cite{DR} for $L(\Fb_a)$ with $a$ possibly negative).
Taking $Q$ to be the hyperfinite II$_1$--factor $R$ and using the result~\cite{D94} that the $n$--fold free product of $R$ is
$L(\Fb_n)$, this gives after rescaling if necessary
that if $1<t<\infty$ and if $\GEu$ is a standard lattice with finite depth, then $L(\Fb_t)$ has
a subfactor whose standard invariant is $\GEu$;
This result was first proved by R\u adulescu~\cite{Ra}.

We also show, (Theorem~\ref{thm:subinf}), that for every standard lattice $\GEu$ and for every II$_1$--factor $Q$,
there is a subfactor $N\subset M$ whose standard lattice is $\GEu$ and such that
\begin{align*}
M&=\staropwl_{i=1}^\infty Q_{s_i} \\[0.5ex]
N&=\staropwl_{i=1}^\infty Q_{t_i}
\end{align*}
for some positive numbers $s_1,s_2,\ldots$ and $t_1,\ldots,t_2$ depending only on $\GEu$.
As a consequence of these two results, we obtain (Theorem~\ref{thm:univ}) that if $Q$ is a II$_1$--factor whose
fundamental group is equal to the set of all positive reals then 
the free product $Q*Q*\cdots$ of infinitely many copies of $Q$ has the universal property with respect to subfactors,
namely, its subfactors give rise to all standard lattices as their standard invariants.

Popa and Shlyakhtenko~\cite{PS} also proved that
if $N\subset M$ is any subfactor of a II$_1$--factor then there is a subfactor pair $\Nh\subseteq\Mh$
having the same standard lattice as $N\subseteq M$ and such that $\Nh\cong N*L(\Fb_\infty)$ and $\Mh\cong M*L(\Fb_\infty)$.
Their proofs involved amalgamated free products of the form
\begin{equation}
\label{eq:IIinffp}
(Q\otimes B,\tau_Q\otimes\id_B)*_B(M\otimes B(H),E)\;,
\end{equation}
where $E$ is a conditional expectation onto a copy of the type I algebra $B$ embedded in the II$_\infty$--factor $M\otimes B(H)$, and where
they use $Q=L(\Fb_\infty)$.

We will identify the free product~\eqref{eq:IIinffp} for an arbitrary II$_1$--factor $Q$ in terms of free subproducts,
and thereby prove
generalizations of the above mentioned result where $L(\Fb_\infty)$ in the formula for $\Nh$ and $\Mh$ is replaced
by free products of rescalings of $Q$.
(See Theorems~\ref{thm:Msubfin}, \ref{thm:Msubinf} and~\ref{thm:NMuniv}.)

The organization of the rest of the paper is as follows:
in~\S\ref{sec:amalgfpI} and, respectively,~\S\ref{sec:amalgfpII}, isomorphism results are proved for amalgamated free products
as in~\eqref{eq:IIinffp} and, respectively~\eqref{eq:MNamalgfp}.
(However, the simplifying assumption is made that the algebra $B$ is commutative --- an easy trick reduces the general case to this.
Moreeover, because it entails no extra difficulty and in fact makes the proof more transparent,
we consider free products where the algebra $Q\otimes B$ of~\eqref{eq:IIinffp} and~\eqref{eq:MNamalgfp} is replaced by a direct sum
$\bigoplus_{i\in I}Q(i)$ of II$_1$--factors whose center is $B$.)
In~\S\ref{sec:subf} we prove the existence of subfactors, including those described above.

\medskip

\noindent
{\em Acknowledgements.}
The author thanks S.~Popa and F.~R\u adulescu for helpful conversations.

\section{Amalgamated free products, I}
\label{sec:amalgfpI}

In this section, we consider amalgamated free products related to those of the form~\eqref{eq:IIinffp} when $B$ is taken to be commutative,
and we show the result can be expressed in terms of free scaled products of von Neumann algebras.
See~\cite{D} for the definition of and results about free scaled products.

Take an index set $I$ equal to $\{1,2,\ldots,m\}$ for some integer $m\ge2$ or to the natural numbers (in which case we write $m=\infty$)
and let $\Qc=\bigoplus_{i\in I}Q(i)$, where each $Q(i)$ is a II$_1$--factor.
let $B$ be the center of $\Qc$ and let $p_i$ denote the minimal projection of $B$ that is the support of the
$i$th summand $Q(i)$ of $\Qc$.
Let $E^\Qc_B:\Qc\to B$ be the conditional expectation that is the center--valued trace on $\Qc$.
Let $M$ be a type II$_1$ or II$_\infty$ factor with fixed faithful (finite or semifinite) trace $\Tr$.
Let $N$ denote the II$_1$--factor obtained as $eMe$ if $M$ is type II$_\infty$, where $e\in M$ is a projection of trace $1$,
and as $M_{\Tr(1)^{-1}}$ if $M$ is type II$_1$;
(thus if $M$ is type II$_1$ and $\Tr(1)\ge1$ then $N\cong eMe$ where $e$ is as above).
Suppose $B$ is embedded in $M$ as a unital subalgebra in such a way that $\Tr(p_i)<\infty$ for all $i\in I$
and let $E^M_B:M\to B$ be the $\Tr$--preserving conditional expectation.
Finally, let $\beta(i)=\Tr(p_i)$.

In this and the next section, when $X_1$ and $X_2$ are subsets of an algebra, we will use the symbols $\Lambdao(X_1,X_2)$ to
denote the set of all alternating words in $X_1$ and $X_2$,
i.e.\ all words of the form $a_1a_2\cdots a_n$, where $a_j\in X_{i(j)}$ and $i(j)\ne i(j+1)$.
In fact, we will frequently blur the distinction between the set of words {\em per se} and the set of elements in the algebra
that are equal to such words.

\begin{prop}
\label{prop:Minffp}
Consider the amalgamated free product of von Neumann algebras
\begin{equation}
\label{eq:inffp}
(\Mcal,E)=(M,E^M_B)*_B(\Qc,E^\Qc_B)\;.
\end{equation}
If $\sum_{i=1}^m\beta(i)<\infty$ then $\Mcal$ is a II$_1$--factor, while if $\sum_{i=1}^m\beta(i)=\infty$ then $\Mcal$ is a II$_\infty$--factor.
In either case, we have
\begin{equation}
\label{eq:p1Mp1}
p_1\Mcal p_1\cong N_{\beta(1)}\staropwl_{i\in I}\big[\tfrac{\beta(i)}{\beta(1)},Q(i)\big]\;.
\end{equation}
\end{prop}
\begin{proof}
We denote also by $\Tr$ the semifinite trace $\Tr\restrict_B\circ E$ on $\Mcal$.
The method of~\cite{D98} shows that $E$ is faithful, hence $\Tr$ is faithful on $\Mcal$;
we will show the isomorphism~\eqref{eq:p1Mp1}, from which it will follow that $\Mcal$ is a factor with faithful
trace $\Tr$, which is finite or semifinite depending on $\sum_{i=1}^m\beta(i)$.

Let $p_{[1,n]}=p_1+p_2+\cdots+p_n$ and let $\beta[1,n]=\Tr(p_{[1,n]})$.
For every $n\in I$, let
\[
\Mcal(n)=W^*(p_{[1,n]}Mp_{[1,n]}\cup p_{[1,n]}\Qc)\;.
\]
Then $\Mcal(1)$ is generated by $p_1Mp_1\cong N_{\beta(1)}$ and $p_1\Qc\cong Q(1)$, which are free, so
\[
\Mcal(1)\cong N_{\beta(1)}*Q(1)\;.
\]

Now for $n\in I$, $n\ge2$, we will find $\Mcal(n)$ in terms of $\Mcal(n-1)$.
Let $K\in\Nats$ be such that $\beta[1,n-1]\ge\beta(n)/K$ and let $(e_{ij})_{1\le i,j\le K}$ be a system of matrix units in $M$ such that $\sum_1^Ke_{ii}=p_n$.
Let $A=\lspan\{e_{ij}\mid 1\le i,j\le K\}$ and let $\Pc=W^*(A\cup p_n\Qc)$.
Clearly $A$ and $p_n\Qc$ are free with respect to $\beta(n)^{-1}\Tr\restrict_{p_n\Mcal p_n}$, so
\[
\Pc\cong\Qc(n)*M_K(\Cpx)\cong\Qc(n)*L(\Fb_{1-K^{-2}}).
\]
Let $v\in M$ be such that $v^*v=e_{11}$ and $q:=vv^*\le p_{[1,n-1]}$\;.
Then 
\[
\Mcal(n)=W^*(\Mcal(n-1)\cup\Pc\cup\{v\})\;, \qquad
q\Mcal(n)q=W^*(q\Mcal(n-1)q\cup v\Pc v^*)\;.
\]

We claim that $q\Mcal(n-1)q$ and $v\Pc v^*$ are free with respect to $K\beta(n)^{-1}\Tr\restrict_{q\Mcal q}$.
This is equivalent to freeness of $v^*\Mcal(n-1)v$ and $e_{11}\Pc e_{11}$
with respect to $K\beta(n)^{-1}\Tr\restrict_{e_{11}\Mcal e_{11}}$, which is what we will show.
Using the trace--preserving conditional expectation $\Pc\to A$ along with the facts that $e_{11}$ is a minimal projection in $A$ and
that $\{p_n\}\cup\Lambdao(A\cap\ker\Tr,p_n\Qc\cap\ker\Tr)$ spans a dense subspace of $\Pc$, we find that
a dense subspace of $e_{11}\Pc e_{11}\cap\ker\Tr$ is spanned by
\[
\Phit:=\bigcup_{1\le i,j\le K}e_{1i}\Phi e_{j1}\;,
\]
where $\Phi$ is the set of all words belonging to
$\Lambdao(A\cap\ker\Tr,p_n\Qc\cap\ker\Tr)$ whose first and last letters come from $p_n\Qc\cap\ker\Tr$.
On the other hand, since
\[
p_{[1,n-1]}B\cup\Lambdao(p_{[1,n-1]}Mp_{[1,n-1]}\cap\ker E,p_{[1,n-1]}\Qc\cap\ker E)
\]
spans a dense subspace of $\Mcal(n-1)$, we see that $v^*\Mcal(n-1)v\cap\ker\Tr$ is contained in the closed linear span of
\[
\Psit:=(e_{11}Me_{11}\cap\ker\Tr)\cup\Psi\;,
\]
where $\Psi$ is the set of all words belonging to
$\Lambdao(M\cap\ker E,p_{[1,n-1]}\Qc\cap\ker E)$ with first letter from $e_{11}Mp_{[1,n-1]}$ and with last letter from
$p_{[1,n-1]}Me_{11}$.
Now to show freeness of $q\Mcal(n-1)q$ and $v\Pc v^*$, it will suffice to show
\begin{equation}
\label{eq:LPP}
\Lambdao(\Phit,\Psit)\subseteq\ker\Tr\;.
\end{equation}
Let $x\in\Lambdao(\Phit,\Psit)$.
If the first letter in $x$ comes from $e_{11}\Phi e_{j1}\subseteq\Phit$ then substitute $e_{11}=(e_{11}-K^{-1}p_n)+K^{-1}p_n$ for this first $e_{11}$,
while if the last letter in $x$ comes from $e_{1i}\Phi e_{11}\subseteq\Phit$ then substitute this value for this last $e_{11}$ and then distribute.
Using $e_{j1}e_{11}Mp_{[1,n-1]}\subseteq\ker E$ and $e_{j1}(e_{11}Me_{11}\cap\ker\Tr)e_{i1}\subseteq\ker E$ for all $i,j\in\{1,\ldots,K\}$
and regrouping, we see that $x$ is equal to a linear combination of at most four words from $\Lambdao(M\cap\ker E,\Qc\cap\ker E)$, and hence $\Tr(x)=0$
by freeness.
We have shown~\eqref{eq:LPP} and thereby freeness of $q\Mcal(n-1)q$ and $v\Pc v^*$.

What we have shown above implies there is an isomorphism
\begin{equation}
\label{eq:Mn1}
p_{[1,n-1]}\Mcal(n)p_{[1,n-1]}\overset\sim\longrightarrow\Mcal(n-1)*\big[\tfrac{\beta(n)}{K\beta[1,n-1]},v\Pc v^*\big]
\end{equation}
intertwining the inclusion $\Mcal(n-1)\hookrightarrow p_{[1,n-1]}\Mcal(n)p_{[1,n-1]}$ and the canonical embedding of $\Mcal(n-1)$ in the RHS of~\eqref{eq:Mn1}.
By results of~\cite{DR},
\[
v\Pc v^*\cong\Pc_{1/K}\cong\Qc(n)_{1/K}*L(\Fb_{K^2-1})\;,
\]
so using the technology of free scaled products~\cite{D},
we get isomorphisms
\begin{equation*}
\begin{aligned}
\Mcal&(n-1)*\big[\tfrac{\beta(n)}{K\beta[1,n-1]},\Qc(n)_{1/K}*L(\Fb_{K^2-1})\big]\overset\sim\longrightarrow \\
&\overset\sim\longrightarrow\Mcal(n-1)*\big[\tfrac{\beta(n)}{K\beta[1,n-1]},\Qc(n)_{1/K}\big]*\big[\tfrac{\beta(n)}{K\beta[1,n-1]},L(\Fb_{K^2-1})\big] \\
&\overset\sim\longrightarrow\Mcal(n-1)*L(\Fb_{(\beta(n)/\beta[1,n-1])^2(1-K^{-2})})*\big[\tfrac{\beta(n)}{K\beta[1,n-1]},\Qc(n)_{1/K}\big] \\
&\overset\sim\longrightarrow\Mcal(n-1)*\big[\tfrac{\beta(n)}{\beta[1,n-1]},\Qc(n)\big]\;,
\end{aligned}
\end{equation*}
all intertwining the canonical embeddings of $\Mcal(n-1)$ into the free scaled products;
(the last isomorphism above is obtained from Theorem~5.5 of~\cite{D}, and the others follow routinely from the definition of free scaled product).
Now compressing by $p_1$ and appealing to Theorem~4.9 of~\cite{D}, we get an isomorphism
\begin{equation}
\label{eq:Mn2}
p_1\Mcal(n)p_1\overset\sim\longrightarrow p_1\Mcal(n-1)p_1*\big[\tfrac{\beta(n)}{\beta(1)},\Qc(n)\big]
\end{equation}
intertwining the inclusion $p_1\Mcal(n-1)p_1\hookrightarrow p_1\Mcal(n)p_1$ and the canonical embedding of $p_1\Mcal(n-1)p_1$ in the RHS of~\eqref{eq:Mn2}.

Taking the inductive limit as $n$ increases gives the desired isomorphism~\eqref{eq:p1Mp1}.
\end{proof}

The following corollary is a direct consequence of Proposition~\ref{prop:Minffp} and results from~\cite{DR} and~\cite{D}.
\begin{cor}
\label{cor:MN}
Let
\[
(\Mcal,E)=(M,E^M_B)*_B(\Qc,E^\Qc_B)
\]
with accompanying notation be as in Proposition~\ref{prop:Minffp}.
\renewcommand{\labelenumi}{(\Alph{enumi})}
\begin{enumerate}

\item
If $I=\{1,\ldots,m\}$ is finite then
\[
\Mcal\cong M*Q(1)_{\frac1{\beta(1)}}*Q(2)_{\frac1{\beta(2)}}*\cdots*Q(m)_{\frac1{\beta(m)}}*L(\Fb_t)
\]
where
\[
t=-m+\sum_{i=1}^m\beta(i)^2\;.
\]

\item
If $I=\Nats$ is infinite and if one or more of the following conditions hold:
\vskip2ex
\renewcommand{\labelenumii}{(\roman{enumii})}
\begin{enumerate}

\item
$Q(i)\cong Q(i)*L(\Fb_\infty)$ for some $i\in\Nats$,

\item
$N\cong N*L(\Fb_\infty)$,

\item
$\displaystyle\sum_{i=1}^\infty\beta(i)^2=\infty$,

\end{enumerate}
\vskip2ex
then
\[
p_1\Mcal p_1\cong N_{\beta(1)}*\bigg(\staropwl_{i=1}^\infty Q(i)_{\frac{\beta(1)}{\beta(i)}}\bigg)\;.
\]

\end{enumerate}
\end{cor}

\section{Amalgamated free products, II}
\label{sec:amalgfpII}

In this section, we will describe amalgamated free product von Neumann algebra related to those in~\eqref{eq:MNamalgfp},
but with $B$ commutative,
in terms of free scaled products.

We begin by letting $A$ be a separable, type I von Neumann algebra with atomic
center $Z(A)$
and let $B$ be a unital commutative subalgebra of $A$ with the property that $B\cap Z(A)=\Cpx1$.
We assume that there is a 
normal, faithful, semifinite trace $\Tr_A$ on $A$ such that $\Tr_A(p)<\infty$ for every minimal projection $p$ in $B$.
Let $E^A_B:A\to B$
be the $\Tr_A$--preserving conditional expectation onto $B$.
Let $(p_i)_{i\in I}$
be the minimal projections in $B$,
where we take either $I=\{1,\ldots,n\}$, some $n\in\Nats$, or $I=\Nats$, in which case we set $n=\infty$.
Let $Q(i)$ be a II$_1$--factor ($i\in I$) and let $\QEu=\bigoplus_{i\in I}Q(i)$,
with the center of $\QEu$ identified with $B$ in such a way
that $p_i$ is the support projection of the $i$th summand $Q(i)$.
Let $E^\QEu_B:\QEu\to B$ be the conditional expectation that is the center valued trace on $\QEu$.
Our goal in this section is to write the amalgamated free product von Neumann algebra
\[
(\PEu,E)=(A,E^A_B)*_B(\QEu,E^\QEu_B)
\]
as a free scaled product of II$_1$--factors.

We have implicitly choosen an ordering of the minimal projections in $B$.
Let $(q_j)_{j\in J}$ be the minimal
projections in $Z(A)$.
Define $J_m$ recursively ($m\in I$) by
\begin{align*}
J_1&=\{j\in J\mid q_jp_1\ne0\} \\
J_m&=\{j\in J\mid q_jp_m\ne0\}\backslash(J_1\cup\cdots\cup J_{m-1})
 \qquad(m\in I\backslash\{1\})\;.
\end{align*}
Let $\beta(i)=\Tr_A(p_i)$, ($i\in I$)
and if $f_j$ is a minimal projection of $A$ that lies under $q_j$, then let
$\alpha(j)=\Tr_A(f_j)$ ($j\in J$).
For every $m\in I\backslash\{1\}$, let $e_m$ be a minimal projection in $A$ 
such that $e_m\le p_1+\cdots+p_{m-1}$ and $e_m$ is equivalent to a subprojection of $p_m$.
(Such a projection must exist by the condition $Z(A)\cap B=\Cpx$.)
Let $j(m)\in J$ be such that $e_m\le q_{j(m)}$ and let $\gamma(m)=\alpha(j(m))=\Tr_A(e_m)$.

\begin{prop}
\label{prop:amalgfp}
Let
\[
(\PEu,E)=(A,E^A_B)*_B(\QEu,E^\QEu_B)
\]
be the free product of von Neumann algebras with amalgamation over $B$.
If $\sum_{i=1}^n\beta(i)<\infty$ then $\PEu$ is a II$_1$--factor, while
if $\sum_{i=1}^n\beta(i)=\infty$ then $\PEu$ is a II$_\infty$--factor.
In either case,
\begin{equation}
\label{eq:p1Miso}
p_1\PEu p_1\cong(Q(1)*L(\Fb_r))
\staropwl_{i=2}^n\,\bigg[\tfrac{\gamma(i)}{\beta(1)},Q(i)_{\frac{\gamma(i)}{\beta(i)}}\bigg]\;,
\end{equation}
where
\[
r=\frac1{\beta(1)^2}\bigg(\beta(1)^2-\sum_{j\in J_1}\alpha(j)^2\bigg)
+\frac1{\beta(1)^2}\sum_{m=2}^n\bigg(\beta(m)^2-\gamma(m)^2-\sum_{j\in J_m}\alpha(j)^2\bigg)\;.
\]
\end{prop}
\begin{proof}
The condition $Z(A)\cap B=\Cpx$ ensures that $\PEu$ is a factor.
We have the normal, faithful, semifinite trace $\Tr=\Tr_A\circ E$ on $\PEu$, which is finite
if and only if $\sum_{i=1}^n\beta(i)<\infty$.

It remains to find the isomorphism class of $p_1\PEu p_1$.
For $m\in I$, let
\begin{align*}
p_{[1,m]}&=p_1+\cdots+p_m \\
\beta[1,m]&=\beta(1)+\cdots+\beta(m)=\Tr(p_{[1,m]})
\end{align*}
and let
\begin{align*}
\NEu(m)&=W^*(p_m \QEu\cup p_mAp_m) \\
\PEu(m)&=W^*(p_{[1,m]}\QEu\cup p_{[1,m]}Ap_{[1,m]})\;.
\end{align*}
Then
\[
\NEu(m)\cong Q(m)*L(\Fb_{s(m)})
\]
where $s(m)$ is the free dimension\footnote{
  We are using the ``free dimension'' that was introduced in~\cite{D93} --- see~\cite{D} for further discusion.
  Since each $p_mAp_m$ is finite dimensional, its free dimension is well defined, and we may use the term ``free dimension''
  without inhibition.}
of $p_mAp_m$ with respect to $\beta(m)^{-1}\Tr\restrict_{p_mAp_m}$.
We have
$\PEu(1)=\NEu(1)$ and
\[
s(1)=\frac1{\beta(1)^2}\bigg(\beta(1)^2-\sum_{j\in J_1}\alpha(j)^2\bigg)\;.
\]

Note that $\PEu$ is the increasing union of $\PEu(1)\subseteq\PEu(2)\subseteq\cdots$ and
thus $p_1\PEu p_1$ is the increasing union of
$p_1\PEu(1)p_1\subseteq p_1\PEu(2)p_1\subseteq\cdots$.
We will show that each $p_1\PEu(m)p_1$ is isomorphic to a
free scaled product of II$_1$--factors in such a way that the
inclusions $p_1\PEu(m)p_1\subseteq p_1\PEu(m+1)p_1$ become canonical embeddings.
Taking the inductive limit will allow us to express $p_1\PEu p_1$ as a free scaled product
of II$_1$--factors.

Fix $m$ such that $m,m+1\in I$.
Let
\begin{align*}
K_m&=\bigcup_{i=1}^m J_i, \\
K_m^{(1)}&=\{k\in K_m\mid q_kp_{m+1}\ne0\} \\
K_m^{(0)}&=K_m\backslash K_m^{(1)}\;.
\end{align*}
The hypothesis $Z(A)\cap B=\Cpx$ ensures that $K_m^{(1)}$ is nonempty;
in fact, $j(m+1)\in K_m^{(1)}$.
Make an ordering $K_m^{(1)}=\{k(i)\mid i\in L\}$ where $L=\{0,1,\ldots,\ell\,\}$ or $L=\{0\}\cup\Nats$,
such that $k(0)=j(m+1)$.
For every $k\in K_m^{(1)}$, let $v_k\in A$ be such that $v_k^*v_k$ is a minimal
projection in $A$, $v_k^*v_k\le p_{[1,m]}$ and $v_kv_k^*\le p_{m+1}$.
Note that $\Tr(v_k^*v_k)=\alpha(k)$;
we assume without loss of generality that $v_{k(0)}^*v_{k(0)}=e_{m+1}$.

Let
\[
\REu(0)=W^*(\PEu(m)\cup\NEu(m+1)\cup\{v_{k(0)}\})\;.
\]
Then
\[
p_{[1,m]}\REu(0)p_{[1,m]}=W^*(\PEu(m)\cup v_{k(0)}^*\NEu(m+1)v_{k(0)})\;.
\]
We will show that the pair
\begin{equation}
\label{eq:MmNm1}
e_{m+1}\PEu(m)e_{m+1},\qquad v_{k(0)}^*\NEu(m+1)v_{k(0)}
\end{equation}
is free with respect to $\gamma(m+1)^{-1}\Tr\restrict_{e_m\PEu e_m}$.
Using $C\oup$ to denote $C\cap\ker E$ for any subalgebra $C$ of $\PEu$,
since $e_m$ is a minimal projection in $A$, $\PEu(m)\oup$ is densely spanned
by
\[
\Theta:=e_m\Big(\Lambdao((\QEu p_{[1,m]})\oup,(p_{[1,m]}Ap_{[1,m]})\oup)
\backslash(p_{[1,m]}Ap_{[1,m]})\oup\Big)e_m
\]
and $\big(v_{k(0)}^*\NEu(m+1)v_{k(0)}\big)\oup$ is densely spanned by
\[
\Omega:=v_{k(0)}^*\Big(\Lambdao((\QEu p_{m+1})\oup,(p_{m+1}Ap_{m+1})\oup)
\backslash(p_{m+1}Ap_{m+1})\oup\Big)v_{k(0)}\;.
\]
To show freeness of the pair~\eqref{eq:MmNm1}, it will suffice to show
$\Lambdao(\Theta,\Omega)\subseteq\ker E$.
However, since
\[
(p_{[1,m]}Ap_{[1,m]})v_{k(0)}^*(p_{m+1}Ap_{m+1})\subseteq\ker E^A_B\;,
\]
given $x\in\Lambdao(\Theta,\Omega)$, by regrouping we see that $x$ is equal to some
$x'\in\Lambdao(\QEu\oup,A\oup)$, which is contained in $\ker E$ by freeness.
This shows that the pair~\eqref{eq:MmNm1} is free.
Therefore, there is an isomorphism
\begin{equation}
\label{eq:P0iso}
p_{[1,m]}\REu(0)p_{[1,m]}\overset\sim\longrightarrow
\PEu(m)*\bigg[\tfrac{\gamma(m+1)}{\beta[1,m]}\,,\,\NEu(m+1)_{\frac{\gamma(m+1)}{\beta(m+1)}}\bigg]
\end{equation}
intertwining the inclusion $\PEu(m)\hookrightarrow p_{[1,m]}\REu(0)p_{[1,m]}$ and
the canonical embedding of $\PEu(m)$ in the RHS of~\eqref{eq:P0iso}.
By~\cite{D}, $\REu(0)$ is a factor.

Define $\REu(i)$ recursively for every $i\in L\backslash\{0\}$ by
\[
\REu(i)=W^*(\REu(i-1)\cup\{v_{k(i)}\}).
\]
We will show that every $\REu(i)$ is a factor that can be written as a free subproduct involving $\REu(i-1)$.
Fixing $i$, suppose that $\REu(i-1)$ is a factor.
Let $w_i\in\REu(i-1)$ be such that $w_i^*w_i=v_{k(i)}^*v_{k(i)}$
and $w_iw_i^*=v_{k(i)}v_{k(i)}^*$.
Let $f=w_i^*w_i$ and $g=w_iw_i^*$;
note that $f\le p_{[1,m]}q_{k(i)}$ and $g\le p_{m+1}q_{k(i)}$.
We will show that in $f\PEu f$ and with respect to $\alpha(k(i))^{-1}\Tr\restrict_{f\PEu f}$,
$w_i^*v_{k(i)}$ is a Haar unitary and is $*$--free from $f\REu(i-1)f$.
Let
\[
C=W^*(p_{[1,m]}Ap_{[1,m]}\cup p_{m+1}Ap_{m+1}\cup\{v_{k(0)},\ldots,v_{k(i-1)}\})\;.
\]
Then $fCg=\{0\}$ and $\REu(i-1)$ is densely spanned by $\Lambdao(C\oup,(p_{[1,m+1]}\QEu)\oup)$.
In order to show that $w_i^*v_{k(i)}$ is a Haar unitary and is $*$--free from $f\REu(i-1)f$,
it will suffice to show
\[
\Lambdao\big((f\REu(i-1)f)\oup\cup f\REu(i-1)g\cup g\REu(i-1)f\cup(g\REu(i-1)g)\oup,
\{v_{k(i)},v_{k(i)}^*\}\big)\subseteq\ker E\;.
\]
Because $f$ and $g$ are minimal projections in $C$ and $fCg=\{0\}$,
\begin{align*}
(f\REu(i-1)f)\oup\cup f\REu(i-1)g\cup g\REu(i-1)f\cup&(g\REu(i-1)g)\oup\subseteq \\
&\subseteq\clspan\big(\Lambdao(C\oup,(p_{[1,m+1]}\QEu)\oup)\backslash C\oup\big)\;.
\end{align*}
But since $Cv_{k(i)}C\subseteq\ker E$, we have
\[
\Lambdao\big(\Lambdao(C\oup,(p_{[1,m+1]}\QEu)\oup)\backslash C\oup,\{v_{k(i)},v_{k(i)}^*\}\big)
\subseteq\ker E\;.
\]
Thus $w_i^*v_{k(i)}$ is a Haar unitary and is $*$--free from $f\REu(i-1)f$.
This gives a $*$--isomorphism
\begin{equation}
\label{eq:Piiso}
p_{[1,m]}\REu(i)p_{[1,m]}\overset\sim\longrightarrow
p_{[1,m]}\REu(i-1)p_{[1,m]}*\bigg[\tfrac{\alpha(k(i))}{\beta[1,m]}\,,\,L(\Ints)\bigg]
\end{equation}
intertwining the inclusion
$p_{[1,m]}\REu(i-1)p_{[1,m]}\hookrightarrow p_{[1,m]}\REu(i)p_{[1,m]}$ and the canonical
embedding of $p_{[1,m]}\REu(i-1)p_{[1,m]}$ into the RHS of~\eqref{eq:Piiso}.

Since $\PEu(m+1)=\overline{\bigcup_{i\in L}\REu(i)}$, by composing the isomorphisms~\eqref{eq:P0iso} and~\eqref{eq:Piiso},
we get an isomorphism
\begin{equation}
\label{eq:Mm1iso}
\begin{aligned}
p_{[1,m]}\PEu(m+1)p_{[1,m]}&\overset\sim\longrightarrow \\[1ex]
&\bigg(\PEu(m)*
\bigg[\tfrac{\gamma(m+1)}{\beta[1,m]}\,,\,\NEu(m+1)_{\frac{\gamma(m+1)}{\beta(m+1)}}\bigg]\bigg)
\staropwl_{i\in L}\,\bigg[\tfrac{\alpha(k(i))}{\beta[1,m]}\,,\,L(\Ints)\bigg]
\end{aligned}
\end{equation}
intertwining the inclusion $\PEu(m)\hookrightarrow p_{[1,m]}\PEu(m+1)p_{[1,m]}$ and the
canonical embedding of $\PEu(m)$ in the RHS of~\eqref{eq:Mm1iso}.
Using that $\NEu(m+1)\cong Q(m+1)*L(\Fb_{s(m+1)})$ with
\[
s(m+1)=1-\sum_{j\in K_m^{(1)}\cup J_m}\bigg(\tfrac{\alpha(j)}{\beta(m+1)}\bigg)^2\;,
\]
by results from~\cite{D} there is an isomorphism
\begin{equation}
\label{eq:p1Piiso}
p_1\PEu(m+1)p_1\overset\sim\longrightarrow
(p_1\PEu(m)p_1*L(\Fb_{r(m+1)}))
*\bigg[\tfrac{\gamma(m+1)}{\beta(1)}\,,\,Q(m+1)_{\frac{\gamma(m+1)}{\beta(m+1)}}\bigg]\;,
\end{equation}
where
\[
r(m+1)=\frac1{\beta(1)^2}\bigg(\beta(m+1)^2-\gamma(m+1)^2-\sum_{j\in J_m}\alpha(j)^2\bigg)\;,
\]
intertwining the inclusion $p_1\PEu(m)p_1\hookrightarrow p_1\PEu(m+1)p_1$ and the canonical
embedding of $p_1\PEu(m)p_1$ in the RHS of~\eqref{eq:p1Piiso}.
Taking the inductive limit as $m\to\infty$, we obtain the isomorphism~\eqref{eq:p1Miso}.
\end{proof}

The next corollary follows directly from Proposition~\ref{prop:amalgfp} and results of~\cite{DR} and~\cite{D}.
\begin{cor}
\label{cor:Mfp}
Let $A$, $B$, $\QEu=\bigoplus_{i\in I}Q(i)$ and
\[
(\PEu,E)=(\QEu,E^\QEu_B)*_B(A,E^A_B)
\]
with accompanying notation be as in Proposition~\ref{prop:amalgfp}.
\renewcommand{\labelenumi}{(\Alph{enumi})}
\begin{enumerate}

\item
If $I=\{1,\ldots,n\}$ is finite then
\[
\PEu\cong Q(1)_{\frac1{\beta(1)}}*Q(2)_{\frac1{\beta(2)}}*\cdots*Q(n)_{\frac1{\beta(n)}}
*L(\Fb_t)
\]
where
\begin{align*}
t&=-n+1+\sum_{i=1}^n\beta(i)^2-\sum_{j\in J}\alpha(j)^2 \\
&=-n+1+\fdim(A)-\fdim(B)\;.
\end{align*}

\item
If $I=\Nats$ is infinite and if one or more of the following conditions hold:
\vskip2ex
\renewcommand{\labelenumii}{(\roman{enumii})}
\begin{enumerate}

\item
$Q(i)\cong Q(i)*L(\Fb_\infty)$ for some $i\in\Nats$,

\item
$\displaystyle\sum_{i=1}^\infty\gamma(i)^2=\infty$,

\item
$r=\infty$,

\end{enumerate}
\vskip2ex
then
\[
p_1\PEu p_1\cong\staropwl_{i=1}^\infty Q(i)_{\frac{\beta(1)}{\beta(i)}}\;.
\]
\end{enumerate}
\end{cor}

\section{Subfactors}
\label{sec:subf}

In this section, we will recall constructions and results of Popa and Shlyakhtenko~\cite{PS}
and we will use them together with the 
results from~\S\ref{sec:amalgfpI} and~\S\ref{sec:amalgfpII} to
prove the existence of subfactors in free subproducts and free products, including those mentioned in the introduction.

Let $Q$ be any II$_1$--factor.
Let $M_{-1}\subset M_0$ be a II$_1$--subfactor of finite index $\lambda$, let $\GEu=\GEu_{M_{-1},M_0}$ be the associated standard $\lambda$--lattice and
let $\Gamma$ be the standard graph and $\Gamma'$ the second standard graph of $M_{-1}\subset M_0$.
\begin{comment}
Let $\sbold$ and $\sbold'$ be the vectors of positive entries so that $\Gamma\Gamma^t\sbold=\lambda^{-1}\sbold$ and
$\Gamma'(\Gamma')^t\sbold'=\lambda^{-1}\sbold'$, normalized so that $\sbold_*=1=\sbold'_*$,
where $*$ denotes the initial vertix of $\Gamma$, respectively of $\Gamma'$.
\end{comment}
Popa and Shlyakhtenko~\cite{PS} construct semifinite von Neumann algebras
\begin{equation}
\label{eq:MAcs}
\begin{matrix}
\MEu_{-1}&\subset&\MEu_0 \\
\cup&&\cup \\
\AEu_{-1}^0&\subset&\AEu_0^0 \\
\cup&&\cup \\
\AEu_{-1}^{-1}&\subset&\AEu_0^{-1}
\end{matrix}
\end{equation}
with a specified faithful semifinite trace $\Tr$ on $\MEu_0$ and they prove the following:
\renewcommand{\labelenumi}{\Roman{enumi}.}
\begin{enumerate}

\item The trace preserving conditional expectations make~\eqref{eq:MAcs} a diagram of commuting squares.

\item
$\MEu_i\cong M_i\otimes B(\HEu)$, ($i=-1,0$).

\item Each algebra $\AEu_i^k$ is a type I von Neumann algebra with atomic center, and $\Tr$ is finite on all minimal projections in $\AEu_i^k$.

\item The graphs of the inclusions $\AEu_{-1}^{-1}\subset\AEu_0^{-1}$ and $\AEu_{-1}^0\subset\AEu_0^0$ are $\Gamma$ and $\Gamma'$, respectively.

\item The commuting square
\begin{equation}
\label{eq:Acs}
\begin{matrix}
\AEu_{-1}^0&\subset&\AEu_0^0 \\
\cup&&\cup \\
\AEu_{-1}^{-1}&\subset&\AEu_0^{-1}
\end{matrix}
\end{equation}
depends only on the standard $\lambda$--lattice $\GEu$, and the commuting square~\eqref{eq:Acs} is functorial in $\GEu$,
(see~\cite[Thm.\ 2.9]{PS} for details).

\item Consider the amalgamated free products of von Neumann algebras
\begin{align}
\PEu_0&=\AEu_0^0*_{\AEu_0^{-1}}(Q\otimes\AEu_0^{-1}) \label{eq:P0} \\
\PEu_{-1}&=\AEu_{-1}^0*_{\AEu_{-1}^{-1}}(Q\otimes\AEu_{-1}^{-1})\;, \label{eq:P-1}
\end{align}
where the amalgamation is with respect to the trace preserving conditional expectations.
Then $\PEu_{-1}\subset\PEu_0$ is an inclusion of type II$_\infty$ factors 
whose standard invariant (i.e.\ the standard invariant
of $p\PEu_{-1}p\subset p\PEu_0p$ where $p\in\PEu_{-1}$ is any finite projection)
is $\GEu$.

\item Consider the amalgamated free products of von Neumann algebras
\begin{align}
\MEuh_0&=\MEu_0*_{\AEu_0^{-1}}(Q\otimes\AEu_0^{-1}) \label{eq:Mh0} \\
\MEuh_{-1}&=\MEu_{-1}*_{\AEu_{-1}^{-1}}(Q\otimes\AEu_{-1}^{-1})\;, \label{eq:Mh-1}
\end{align}
where the amalgamation is with respect to the trace preserving conditional expectations.
Then $\MEuh_{-1}\subset\MEuh_0$ is an inclusion of type II$_\infty$ factors 
whose standard invariant (i.e.\ the standard invariant
of $p\MEuh_{-1}p\subset p\MEuh_0p$ where $p\in\MEuh_{-1}$ is any finite projection)
is $\GEu$.

\end{enumerate}

\begin{numbered}\rm
\label{num:P-1P0}
Consider the commuting square~\eqref{eq:Acs} determined by $\GEu$.
Let $\PEu_0$ be as in~\eqref{eq:P0} and let $q(0)\in\Ac_0^{-1}$ be the sum of a maximal family of mutually orthogonal and nonequivalent
minimal projections of $\Ac_0^{-1}$.
Letting $B=q(0)\Ac_0^{-1}q(0)$, we have that $B$ is commutative and
\begin{equation*}
q(0)\PEu_0q(0)=q(0)\Ac_0^0q(0)*_B(Q\otimes B)\;,
\end{equation*}
where the amalgamation is with respect to the trace preserving conditional expectations.
Fix any minimal projection $q_1(0)$ in $B$.
Let $m(0)\in\Nats\cup\{\infty\}$ be the number of minimal projections in $B$
and let $(\beta(0,i))_{i=1}^{m(0)}$ be the values of $\Tr$ on the minimal projections in $B$, with $\beta(0,1)=\Tr(q_1(0))$.
Proposition~\ref{prop:amalgfp} shows that for cretain $r(0)\in[0,\infty]$ and $\gamma(0,i)\le\beta(0,i)$, we have
\begin{equation}
\label{eq:qP0q}
q_1(0)\PEu_0q_1(0)\cong\begin{cases}
Q*L(\Fb_{r(0)})&\text{if }m(0)=1 \\
\big(Q*L(\Fb_{r(0)})\big)\displaystyle\staropwl_{i=2}^{m(0)}\bigg[\frac{\gamma(0,i)}{\beta(0,1)},Q_{\frac{\gamma(0,i)}{\beta(0,i)}}\bigg]&\text{if }m(0)>1\;.
\end{cases}
\end{equation}
Analogously, if $\PEu_{-1}$ is as in~\eqref{eq:P-1} and 
if $q_1(-1)$ is any minimal projection in $\Ac_{-1}^{-1}$, if $m(-1)$ is the 
number of minimal projections in the center of $\Ac_{-1}^{-1}$, letting $(\beta(-1,i))_{i=1}^{m(-1)}$
be $\Tr$ applied to a maximal family of mutually inequivalent minimal projections in $\Ac_{-1}^{-1}$, with $\beta(-1,i)=\Tr(q_1(-1))$,
we have
\begin{equation}
\label{eq:qP-1q}
q_1(-1)\PEu_{-1}q_1(-1)\cong\begin{cases}
Q*L(\Fb_{r(-1)})&\text{if }m(-1)=1 \\
\big(Q*L(\Fb_{r(-1)})\big)\displaystyle\staropwl_{i=2}^{m(-1)}\bigg[\frac{\gamma(-1,i)}{\beta(-1,1)},Q_{\frac{\gamma(-1,i)}{\beta(-1,i)}}\bigg]&\text{if }m(-1)>1
\end{cases}
\end{equation}
for certain $r(-1)\in[0,\infty]$ and $\gamma(-1,i)\le\beta(-1,i)$.
\end{numbered}

\begin{thm}
\label{thm:subfin}
Let $\GEu$ be a standard lattice of finite depth.
Then there are $m,n\in\Nats$, positive real numbers $s_1,\ldots,s_m$ and $t_1,\ldots,t_n$
and there are $a>1-m$ and $b>1-n$ such that for any II$_1$--factor $Q$, there is an inclusion $P_{-1}\subset P_0$ of II$_1$--factors
whose standard invariant is $\GEu$ and such that
\begin{align*}
P_0&=Q_{s_1}*\cdots*Q_{s_m}*L(\Fb_a) \\
P_{-1}&=Q_{t_1}*\cdots*Q_{t_n}*L(\Fb_b)\;.
\end{align*}
\end{thm}
\begin{proof}
Since $\GEu$ is of finite depth, it follows from~\cite{PS} (see point IV above) that the centers of $\Ac_0^{-1}$ and $\Ac_{-1}^{-1}$
are finite dimensional.
Thus $m(0)$ and $m(1)$ in~\eqref{eq:qP0q} and~\eqref{eq:qP-1q} are finite.
By results in~\cite{D}, from~\eqref{eq:qP0q} and~\eqref{eq:qP-1q} we get
\begin{align*}
q_1(0)\PEu_0q_1(0)&\cong Q_{s_1'}*\cdots*Q_{s_{m(0)}'}*L(\Fb_{a'}) \\
q_1(-1)\PEu_{-1}q_1(-1)&\cong Q_{t_1'}*\cdots*Q_{t_{m(1)}'}*L(\Fb_{b'})
\end{align*}
for some $a'>1-m(0)$ and $b'>1-m(-1)$,
where $s_i'=\beta(0,i)/\beta(0,1)$ and $t_i'=\beta(-1,i)/\beta(-1,1)$.
By results of~\cite{PS} (see point VI above), the II$_1$--subfactor
$q\PEu_{-1}q\subset q\PEu_0q$
has standard invariant $\GEu$, where $q$ is any finite projection in $\PEu_{-1}$.
But $q\PEu_0q$ is a rescaling of $q_1(0)\PEu_0q_1(0)$ by $\lambda:=\Tr(q)/\Tr(q_1(0))$,
and thus by results of~\cite{DR} we have
\[
q\PEu_0q\cong Q_{\lambda s_1'}*\cdots*Q_{\lambda s_{m(0)}'}*L(\Fb_a)
\]
where $a=\lambda^{-2}a'+(m(0)-1)(\lambda^{-2}-1)$.
Rescaling $q_1(-1)\PEu_{-1}q_1(-1)$ yields
\[
q\PEu_{-1}q\cong Q_{t_1}*\cdots*Q_{t_{m(1)}}*L(\Fb_b)
\]
for appropriate values of $t_i$ and $b$.
\end{proof}

\begin{remark}\rm
In the proof of Theorem~\ref{thm:subfin}, $\lambda$ can be chosen to arrange either $a=0$, provided $m(0)>1$,
or $b=0$, provided $m(1)>1$.
\end{remark}

\begin{thm}
\label{thm:subinf}
Given a standard lattice $\GEu$, there are some positive real numbers $s_1,s_2,\ldots$ and $t_1,t_2,\ldots$
such that for any II$_1$--factor $Q$, there is an inclusion $P_{-1}\subset P_0$ of II$_1$--factors
whose standard invariant is $\GEu$ and such that
\begin{equation}
\label{eq:P0P-1}
\begin{aligned}
P_0&\cong\staropwl_{i=1}^\infty Q_{s_i} \\[0.5ex]
P_{-1}&\cong\staropwl_{i=1}^\infty Q_{t_i}
\end{aligned}
\end{equation}
\end{thm}
\begin{proof}
If $\GEu$ is of finite depth, let $\Qt=Q*Q*\cdots$ be a the free products of infinitely many copies of $Q$.
By~\cite{DR00}, a free product of infinitely many II$_1$--factors absorbs a free product with $L(\Fb_\infty)$,
and the rescaling by $\lambda$ of a free product of infinitely many factors is isomorphic to the free product of the same factors
each rescaled by $\lambda$.
Applying these facts and Theorem~\ref{thm:subfin} to $\Qt$,
we find a subfactor $P_{-1}\subset P_0$ whose standard invariant is $\GEu$
and where
\begin{align*}
P_0&\cong(Q_{s_1}*Q_{s_1}*\cdots)*\cdots*(Q_{s_n}*Q_{s_n}*\cdots) \\
P_{-1}&\cong(Q_{t_1}*Q_{t_1}*\cdots)*\cdots*(Q_{t_m}*Q_{t_m}*\cdots)\;.
\end{align*}

If $\GEu$ has infinite depth, let $\Qt=Q*L(\Fb_\infty)$.
Let $\PEu_0$ and $\PEu_{-1}$ be as in~\eqref{eq:P0} and~\eqref{eq:P-1}, but with $Q$ replaced by $\Qt$.
Then we have the isomorphisms as in~\eqref{eq:qP0q} and~\eqref{eq:qP-1q} but with $Q$ replaced by $\Qt$.
Since $\GEu$ has infinite depth, $m(0)$ and $m(-1)$ are infinite.
Since $\Qt$ absorbs a free product with $L(\Fb_\infty)$, using~\cite[Theorem 3.11]{D}
and the above mentioned results about rescalings and absorbtion of $L(\Fb_\infty)$, we have the following isomorphisms:
\[
\begin{matrix}
q_1(0)\PEu_0q_1(0)&\cong&\Qt*\bigg(\displaystyle\staropwl_{i=2}^\infty\Qt_{\frac{\beta(0,i)}{\beta(0,1)}}\bigg)
 &\cong&\displaystyle\staropwl_{i=1}^\infty Q_{\frac{\beta(0,i)}{\beta(0,1)}} \\[3ex]
q_1(-1)\PEu_{-1}q_1(-1)&\cong&\Qt*\bigg(\displaystyle\staropwl_{i=2}^\infty\Qt_{\frac{\beta(-1,i)}{\beta(-1,1)}}\bigg)
 &\cong&\displaystyle\staropwl_{i=1}^\infty Q_{\frac{\beta(-1,i)}{\beta(-1,1)}}
\end{matrix}
\]
As in the proof of Theorem~\ref{thm:subfin}, we get a subfactor $P_{-1}=q\PEu_{-1}q\subset P_0=q\PEu_0q$ with isomorphisms~\eqref{eq:P0P-1}
as desired.
\end{proof}

The following theorem shows that certain II$_1$--factors 
other than $L(\Fb_\infty)$ have the universal property for subfactors.

\begin{thm}
\label{thm:univ}
Let $Q$ be any II$_1$--factor whose fundamental group is equal to the positive real numbers,
and let $P=Q*Q*\cdots$ be the free product of infinitely many copies of $Q$.
Then for any standard lattice $\GEu$ there is a subfactor $N=N_\GEu$ of $P$ whose 
standard invariant is $\GEu$;
also $N$ is isomorphic to $P$.
Moreover, the map $\GEu\mapsto N_\GEu\subset P$ is functorial, as described in~\cite[Theorem 4.3]{PS}.
\end{thm}
\begin{proof}
The existence of the subfactor $N_\GEu\subset P$ is a direct application of Theorem~\ref{thm:subinf}.
Functoriality follows from results of Popa and Shlyakhtenko as in the proof of~\cite[Theorem 4.3]{PS}.
\end{proof}

Consider the commuting square~\eqref{eq:MAcs} determined by the subfactor $M_{-1}\subset M_0$
and let $\MEuh_0$ be as in~\eqref{eq:Mh0}.
With notation as in~\ref{num:P-1P0}, we have
\begin{align*}
q(0)\MEuh_0q(0)&=q(0)\MEu_0q(0)*_B(Q\otimes B), \\
q(-1)\MEuh_{-1}q(-1)&=q(-1)\MEu_{-1}q(-1)*_B(Q\otimes B)\;,
\end{align*}
where the amalgamated free products are with respect to $\Tr$--preserving conditional expectations.
From Proposition~\ref{prop:Minffp}, we have
\begin{equation*}
q_1(0)\MEuh_0q_1(0)\cong\begin{cases}
(M_0)_{\beta(0,1)}&\text{if }m(0)=1 \\
(M_0)_{\beta(0,1)}\displaystyle\staropwl_{i=2}^{m(0)}\big[\tfrac{\beta(0,i)}{\beta(0,1)},Q\big]&\text{if }m(0)>1
\end{cases}
\end{equation*}
and similarly
\begin{equation*}
q_1(-1)\MEuh_0q_1(-1)\cong\begin{cases}
(M_{-1})_{\beta(-1,1)}&\text{if }m(-1)=1 \\
(M_{-1})_{\beta(-1,1)}\displaystyle\staropwl_{i=2}^{m(-1)}\big[\tfrac{\beta(-1,i)}{\beta(-1,1)},Q\big]&\text{if }m(-1)>1\;.
\end{cases}
\end{equation*}
From these, and in a manner analogous to the proofs of Theorems~\ref{thm:subfin}, \ref{thm:subinf} and~\ref{thm:univ},
one proves the following results.

\begin{thm}
\label{thm:Msubfin}
Let $Q$ be any II$_1$--factor.
Then for any finite depth subfactor $M_{-1}\subset M_0$ of II$_1$--factors, there are $m,n\in\Nats$, positive real
numbers $s_1,\ldots,s_m$ and $t_1,\ldots,t_n$ and there are $a>1-m$ and $b>1-n$, all depending only on the standard invariant
of $M_{-1}\subset M_0$,
and there is a II$_1$--subfactor $\Mh_{-1}\subset\Mh_0$
having the same standard invariant as $M_{-1}\subset M_0$ and where
\begin{align*}
\Mh_0&\cong M_0*Q_{s_1}*\cdots*Q_{s_m}*L(\Fb_a) \\
\Mh_{-1}&=M_{-1}*Q_{t_1}*\cdots*Q_{t_n}*L(\Fb_b)\;.
\end{align*}
\end{thm}

\begin{thm}
\label{thm:Msubinf}
Let $Q$ be any II$_1$--factor.
Then for any finite depth subfactor $M_{-1}\subset M_0$ of II$_1$--factors, there are positive real
numbers $s_1,\ldots,s_m$ and $t_1,\ldots,t_n$ depending only on the standard invariant
of $M_{-1}\subset M_0$
and there is a II$_1$--subfactor $\Mh_{-1}\subset\Mh_0$
having the same standard invariant as $M_{-1}\subset M_0$ and where
\begin{align*}
\Mh_0&\cong M_0*\bigg(\staropwl_{i=1}^\infty Q_{s_i}\bigg) \\
\Mh_{-1}&=M_{-1}*\bigg(\staropwl_{i=1}^\infty Q_{t_i}\bigg)\;.
\end{align*}
\end{thm}

\begin{thm}
\label{thm:NMuniv}
Let $Q$ be any II$_1$--factor having fundamental group equal to $\Reals_+^*$
and let $P=Q*Q*\cdots$ be the free product of infinitely many copies of $Q$.
Then for any subfactor $N\subset M$ of a II$_1$--factor,
there is a subfactor $\Nh\subset\Mh$ whose standard lattice is equal to that of $N\subset M$ and such that
$\Nh\cong N*P$ and $\Mh\cong M*P$.
\end{thm}

% \newpage

\bibliographystyle{plain}

\begin{thebibliography}{9}

\bibitem{D93} K.\ Dykema,
{\em Free products of hyperfinite von Neumann algebras and free dimension,}
Duke Math.\ J.\ {\bf 69} (1993), 97-119.

\bibitem{D94} \rule{3em}{.1mm},
{\em Interpolated free group factors,}
Pacific J.\ Math.\ {\bf 163} (1994), 123-135.

\bibitem{D98} \rule{3em}{.1mm},
{\em Faithfulness of free product states,}
J.\ Funct.\ Anal.\ {\bf 154} (1998), 223-229.

\bibitem{D} \rule{3em}{.1mm},
{\em Free subproducts and free scaled products of II$_1$--factors,}
J.\ Funct.\ Anal.\ (to appear).

\bibitem{DR00} K.\ Dykema, F.\ R\u adulescu,
{\em Compressions of free products of von Neumann algebras,}
Math.\ Ann.\ {\bf 316} (2000), 61-82.

\bibitem{DR} \rule{3em}{.1mm},
{\em Rescalings of free products of II$_1$--factors,}
Proc.\ Amer.\ Math.\ Soc.\ (to appear).

\bibitem{P} S.\ Popa,
{\em An axiomatization of the lattice of higher relative commutants of a subfactor,}
Invent.\ Math.\ {\bf 120} (1995), 427-445.

\bibitem{PS} S.\ Popa, D.\ Shlyakhtenko,
{\em Universal properties of $L(\Fb_\infty)$ in subfactor theory,}
preprint (2000).

\bibitem{Ra} F.\ R\u adulescu,
{\em Random matrices, amalgamated free products and subfactors of the von Neumann algebra of a free group, of noninteger index,}
Invent.\ Math.\ {\bf 115} (1994), 347-389.

\end{thebibliography}

\end{document}